\documentclass[11pt]{amsart}

\newcommand{\eps}{\epsilon}
\newcommand{\pa}{\partial}

\newfont{\fnt}{cmr10 scaled 550}
\renewcommand{\eps}{\varepsilon} 

\newtheorem{theorem}{Theorem}[section]

\newtheorem{conj}{Conjecture}
\newtheorem{lemma}{Lemma}[section]

\newtheorem{prop}{Proposition}[section]

\theoremstyle{remark}
\newtheorem{rem}{Remark}[section] 
\numberwithin{equation}{section}

\usepackage{pdfsync}
\font\strange=msbm10

\renewcommand{\epsilon}{\varepsilon}

\renewcommand{\Sigma}{\varSigma}
\newcommand{\R}{{{\mathchoice  {\hbox{$\textstyle{\text{\strange R}}$}}
{\hbox{$\textstyle{\text{\strange R}}$}}
{\hbox{$\scriptstyle  N\kern-0.3em  R$}}  
{\hbox{$\scriptscriptstyle  R\kern-0.2em  R$}}}}}

\newcommand{\Z}{{{\mathchoice  {\hbox{$\textstyle{\text{\strange Z}}$}}
{\hbox{$\textstyle{\text{\strange Z}}$}}
{\hbox{$\scriptstyle  Z\kern-0.3em  Z$}}
{\hbox{$\scriptscriptstyle  Z\kern-0.2em  Z$}}}}}

\newcommand{\N}{{{\mathchoice  {\hbox{$\textstyle{\text{\strange N}}$}}
{\hbox{$\textstyle{\text{\strange N}}$}}
{\hbox{$\scriptstyle  N\kern-0.3em  N$}}
{\hbox{$\scriptscriptstyle  N\kern-0.2em  N$}}}}}

\renewcommand{\phi}{\varphi}
\usepackage{amsmath,amsthm,amscd}

\begin{document}
\title[Spectrum of quantum layers]{On the ground state of quantum layers}
\date{July 21, 2007}
\thanks{The
 author is partially supported by  NSF Career award DMS-0347033 and the
Alfred P. Sloan Research Fellowship.}

 \author{Zhiqin Lu}

\address{Department of
Mathematics, University of California,
Irvine, Irvine, CA 92697}

\email[Zhiqin Lu]{zlu@math.uci.edu}

\maketitle

\pagestyle{myheadings}

\newcommand{\nn}{{\bf n}}
\newcommand{\ka}{K\"ahler }
\newcommand{\ii}{\sqrt{-1}}

\section{Introduction}
The problem is from mesoscopic physics: let $p: \Sigma\to\mathbb R^3$ be an embedded surface in $\mathbb R^3$, we assume that
\begin{enumerate}
\item $\Sigma$ is orientable, complete, but non-compact;
\item $\Sigma$ is not totally geodesic;
\item $\Sigma$ is asymptotically flat in the sense that the second fundamental form goes to zero at infinity.
\end{enumerate}

On can build  a quantum layer $\Omega$ over such a surface $\Sigma$ as follows: as a differentiable manifold, $\Omega=\Sigma\times [-a,a]$ for some positive number $a$. Let $\vec N$ be the unit normal vector of $\Sigma$ in $\mathbb R^3$. Define
\[
\tilde p: \Omega\to\mathbb R^3
\]
by 
\[
\tilde p(x,t)=p(x)+t\vec N_x.
\]
Obviously, if $a$ is small, then $\tilde p$ is an embedding. The Riemannian metric $ds_\Omega^2$ is defined as the pull-back of the Euclidean metric via $\tilde p$. The Riemannian manifold $(\Omega, ds^2_\Omega)$ is called the quantum layer.

Let $\Delta=\Delta_\Omega$ be the Dirichlet Laplacian. Then we make the following

\begin{conj}\label{ccon}
Using the above notations, and further assume that
\begin{equation}\label{huber}
\int_\Sigma|K|d\Sigma<+\infty.
\end{equation}
Then the ground state of $\Delta$ exists.
\end{conj}

We make the following explanation of the notations and terminology:

\begin{enumerate}
\item $\Omega$ is a smooth manifold with boundary. The Dirichlet Laplacian is the self-adjoint extenstion of the Laplacian acting on $C_0^\infty(\Omega)$;
\item By a theorem of Huber~\cite{huber}, if ~\eqref{huber} is valid, then $\Sigma$ is differmorphic to a compact Riemann surface with finitely many points removed. Moreover, White~\cite{bwhite}  proved that if 
\[
\int_\Sigma K^-d\Sigma<+\infty,
\]
then 
\[
\int_\Sigma|K|d\Sigma<+\infty.
\]
Thus~\eqref{huber} can be weakened.
\item Since $\Delta$ is a self-adjoint operator, the spectrum of $\Delta$ is the disjoint union of two parts: pure point spectrum (eigenvalues of finite multiplicity) and the essential spectrum. The ground state is the smallest eigenvalue with finite multiplicity.
\item The conjecture was proved under the condition
\[
\int_\Sigma Kd\Sigma\leq 0
\]
in~\cite{dek,cek-1} by Duclos, Exner and Krej\v{c}i\v{r}{\'\i}k and later by Carron, Exner, and Krej\v{c}i\v{r}{\'\i}k. Thus the remaining case is when
\[
\int_\Sigma Kd\Sigma>0.
\]
By a theorem of Hartman~\cite{hartman}, we know that
\[
\int_\Sigma K=2\pi\chi(\Sigma)-\sum\lambda_i
\]
where $\lambda_i$ are the isoperimetric constants at each end of $\Sigma$. Thus we have
\[
\chi(\Sigma)>0
\]
and $g=0$. The surface must be differmorphic to $\mathbb R^2$. However, even through the topology of the surface is completely known, this is the most difficult case for the conjecture.
\end{enumerate}

\section{Variational Principle}
It is well known that
\[
\sigma_0=\underset{f\in C_0^\infty(\Omega)}{\inf}\frac{\int_\Omega|\nabla f|^2 d\Omega}{\int_\Omega f^2 d\Omega}
\]
is the infimum of the Laplacian, and 
\begin{equation}\label{ess}
\sigma_{ess}=\underset{K}{\sup}\underset{f\in C_0^\infty(\Omega\backslash K)}{\inf}\frac{\int_\Omega|\nabla f|^2 d\Omega}{\int_\Omega f^2 d\Omega}
\end{equation}
is the infimum of the essential spectrum,
where $K$ is running over all the compact subset of $\Omega$. 
Since $\Omega=\Sigma\times [-a,a]$, it is not hard to see that
\begin{equation}\label{ess-1}
\sigma_{ess}=\underset{K\subset\Sigma}{\sup}\,\,\underset{f\in C_0^\infty(\Omega\backslash K\times[-a,a])}{\inf}\frac{\int_\Omega|\nabla f|^2 d\Omega}{\int_\Omega f^2 d\Omega},
\end{equation}
where $K$ is running over all the compact set of $\Sigma$.

By  definition, we have $\sigma\leq\sigma_{ess}$. Furthermore, we have

\begin{prop}
If $\sigma_0<\sigma_{ess}$, then the ground state exists and is equal to $\sigma_0$.
\end{prop}

\qed

Let $(x_1,x_2)$ be a local coordinate system of $\Sigma$. Then $(x_1,x_2,t)$ defines a local coordinate system of $\Omega$. Such a local coordinate system is called a {\it Fermi coordinate system}. Let $x_3=t$ and let $ds_\Omega^2=G_{ij}dx_idx_j$. Then  we have
\begin{equation}\label{gij}
G_{ij}=\left\{
\begin{array}{ll}
(p+t\vec N)_{x_i}(p+t\vec N)_{x_j} & 1\leq i,j\leq 2;\\
0& i=3, \text {or } j=3, \text{ but } i\neq j;\\
1& i=j=3.
\end{array}
\right.
\end{equation}

We make the following defintion: let $f$ be a smooth function of $\Omega$. Then we define
\begin{align}
&Q(f,f)=\int_\Omega|\nabla f|^2 d\Omega-\kappa^2\int_{\Omega} f^2 d\Omega;\label{q}\\
& Q_1(f,f)=\int_\Omega|\nabla' f|^2 d\Omega;\label{q1}\\
& Q_2(f,f)=\int_\Omega\left(\frac{\pa f}{\pa t}\right)^2 d\Omega-\kappa^2\int_{\Omega} f^2 d\Omega,\label{q2}
\end{align}
where $|\nabla' f|^2=\sum_{i,j=1}^2 G^{ij}\frac{\pa f}{\pa x_i}\frac{\pa f}{\pa x_j}$.

Obviously, we have
\[
Q(f,f)=Q_1(f,f)+Q_2(f,f).
\]

It follows that
\[
\int_\Omega|\nabla f|^2 d\Omega=\int_\Omega|\nabla' f|^2 d\Omega
+\int_\Omega\left(\frac{\pa f}{\pa t}\right)^2 d\Omega
\]
for a smooth function $f\in C^\infty(\Omega)$, where
\[
|\nabla' f|^2=\sum_{1\leq i,j\leq 2} G^{ij}\frac{\pa f}{\pa x_i}\frac{\pa f}{\pa x_j}
\]
is the norm of the horizontal differential. Apprently, we have
\[
\int_\Omega|\nabla f|^2 d\Omega\geq
\int_\Omega\left(\frac{\pa f}{\pa t}\right)^2 d\Omega. 
\]

Let  $ds^2_\Sigma=g_{ij} dx_idx_j$  be the Riemannian metric of $\Sigma$ under the coordinates $(x_1,x_2)$. Then we are above to compare the matrices $(G_{ij})_{1\leq i,j\leq 2}$ and $(g_{ij})$, at least outside a big compact set of $\Sigma$.
By~\eqref{gij}, we have
\[
G_{ij}=g_{ij}+tp_{x_i}\vec N_{x_j}+tp_{x_j}\vec N_{x_i}+t^2\vec N_{x_i}\vec N_{x_j}.
\]
We assume that at the point $x$, $g_{ij}=\delta_{ij}$. Then we have
\[
|G_{ij}-\delta_{ij}|\leq 2a|B|+a^2|B|^2,
\]
where $B$ is the second fundamental form of the surface $\Sigma$. Thus we have the following conclusion:

\begin{prop}\label{prop22}
 For any $\eps>0$, there is a compact set $K$ of $\Sigma$ such that on $\Sigma\backslash K$ we have
\[
(1-\eps)
\begin{pmatrix}
g_{11}&g_{12}\\
g_{21}& g_{22}
\end{pmatrix}
\leq
\begin{pmatrix}
G_{11}&G_{12}\\
G_{21}& G_{22}
\end{pmatrix}
\leq (1+\eps)
\begin{pmatrix}
g_{11}&g_{12}\\
g_{21}& g_{22}
\end{pmatrix}.
\]
In particular, we have
\[
(1-\eps)^2d\Sigma dt\leq d\Omega\leq(1+\eps)^2 d\Sigma dt.
\]
\end{prop}

\qed

Let $\kappa=\frac{\pi}{2a}$. Then we proved the followng:

\begin{lemma}\label{lemess}
 Using the above notations, we have
\[
\sigma_{ess}\geq \frac{\pi^2}{4a^2}.
\]
\end{lemma}

{\bf Proof.}  Let $K$ be any compact set of $\Sigma$. If $f\in C_0^\infty(\Omega\backslash K)$, then by Proposition~\ref{prop22}, we have
\[
\int_{\Omega}\left(\frac{\pa f}{\pa t}\right)^2 d\Omega\geq (1-\eps)^2
\int_\Sigma\int_{-a}^a \left(\frac{\pa f}{\pa t}\right)^2 dtd\Sigma
\geq (1-\eps)^2\kappa^2
\int_\Sigma\int_{-a}^a f^2 dtd\Sigma,
\]
where the last inequality is from the $1$-dimensional Poincar\'e inequality.
Thus by using Proposition~\ref{prop22} again, we have
\[
\int_\Omega |\nabla f|^2 d\Omega\geq\frac{(1-\eps)^2}{(1+\eps)^2}
\kappa^2\int_\Omega f^2 dtd\Sigma.
\]
for any $\eps$. 
Thus we have
\[
\sigma_{ess}\geq \frac{(1-\eps)^2}{(1+\eps)^2}\kappa^2
\]
 and the lemma is proved.
 
 \qed
 
 \begin{rem}\label{rem21} Although not needed in this paper, we can actually prove that 
 $\sigma_{ess}=\kappa^2$. To see this, we first observe that  since the second fundamental form of $\Sigma$ is bounded, there is a lower bound for the injectivity radius. As a  result, the volume of the surface $\Sigma$ is infinite. By the assumption, the Gauss curvature is integrable. Thus $\Sigma$ is {\it parabolic} (cf.~\cite{li}). From the above, we conclude that for any $\eps,C>0$ and any compact sets $K\subset\subset K'$ of $\Sigma$, there is a smooth function $\phi\in C_0^\infty(\Sigma\backslash K')$ such that
 \[
 \phi\equiv 1\text{ on } K,\qquad \int_\Sigma\phi^2 d\Sigma>C,\qquad\text {and }
 \int_\Sigma |\nabla\phi|^2 d\Sigma<\eps.
 \]

 Let $\tilde\phi=\phi\chi$, where $\chi=\cos\kappa t$. Then $\tilde\phi$ is a function on $\Omega$ with compact support. Since the second fundamental form goes to zero at infinity, by Proposition~\ref{prop22}, for $K'$ large enough, we have
 \[
 \int_\Omega|\nabla\tilde\phi|^2 d\Omega<4a\eps.
 \]
 Thus from~\eqref{q} and Proposition~\ref{prop22} again, we have
 \begin{align*}
& Q(\tilde\phi,\tilde\phi)<4a\eps+(1+\eps)^2\int_\Sigma\phi^2 d\Sigma\int^a_{-a}\left(\frac{\pa\chi}{\pa t}\right)^2 dt\\&\qquad-(1-\eps)^2\kappa^2\int_\Sigma\phi^2d\Sigma\int^a_{-a}\chi^2 dt.
 \end{align*}
 A straightforward computation gives
 \[
 \int^a_{-a}\left(\frac{\pa\chi}{\pa t}\right)^2 dt=\kappa^2\int^a_{-a}\chi^2 dt.
 \]
 Thus
 \[
 Q(\tilde\phi,\tilde\phi)\leq 4a\eps+4a\eps\int_\Sigma\phi^2d \Sigma.
 \]
 By the definition of $\sigma_{ess}$, we have
 \[
 \sigma_{ess}-\kappa^2\leq\frac{Q(\tilde\phi,\tilde\phi)}{\int_\Omega\tilde\phi^2 d\Omega}\leq\frac{4\eps(1+\int_\Sigma\phi^2 d\Sigma)}{(1-\eps)^2\int_\Sigma\phi^2 d\Sigma}.
 \]
 We let $\eps\to 0$ and $C\to\infty$,  then we have $\sigma_{ess}\leq\kappa^2$, as needed.
  \end{rem}
 
\section{The upper bound of $\sigma_0$}

It is usually more difficult to estimate $\sigma_0$ from above. In ~\cite[Theorem 1.1]{ll-1}, we proved the following
\begin{theorem}\label{main3}
Let $\Sigma$ be a convex surface in $\R^3$ which can be
represented by the graph of a  convex function
$z=f(x,y)$. Suppose $0$ is the minimum point of
the function and suppose that at $0$, $f$ is
strictly convex. Furthermore suppose that the
second fundamental form goes to zero at infinity.
Let $C$ be the supremum of the second fundamental
form of $\Sigma$. Let $Ca<1$. Then the ground
state of the quantum layer $\Omega$  exists.
\end{theorem}

In this section, we generalize the above result into the following:

\begin{theorem}\label{main4}
Let $\Sigma$ be a complete surface in $\R^3$ with nonnegative Gauss curvature but not totally geodesic. Furthermore suppose that the
second fundamental form of $\Sigma$ goes to zero at infinity.
Let $C$ be the supremum of the second fundamental
form of $\Sigma$. Let $Ca<1$. Then the ground
state of the quantum layer $\Omega$  built over $\Sigma$ with width $a$ exists.
\end{theorem}

\begin{rem} Since for all convex function $f$ in Theorem~\ref{main3}, the Gauss curvature is nonnegative, the above theorem is indeed a generalization of Theorem~\ref{main3}. On the other hand, by a theorem of Sacksteder~\cite{sack}, any complete surface of nonnegative curvature is either a developable surface or the graph of some convex function. At a first glance, it seems that there is not much difference between the surfaces in both theorems. However, we have to use a complete different method to prove this slight generalization.
\end{rem}
 
 {\bf Proof of Theorem~\ref{main4}.} 
 If the Gauss curvature is identically zero, then By ~\cite[Theorem 2]{ll-4}, the ground state exists.
 
 If the Guass curvature is positive at one point, then 
 by using the theorem of Sacksteder~\cite{sack}, $\Sigma$ can be represented by the graph of some convex function. If we fix an orientation, we can assume that $H$, the mean curvature, is always nonnegative.
 
 By a result of White~\cite{bwhite}, we know that there is an $\eps_0>0$
 such that for $R>>0$,
 \[
 \int_{\pa B(R)}||B||>\eps_0,
 \]
 where $B$ is the second fundamental form of $\Sigma$. Since $\Sigma$ is convex, we have
 \[
 H\geq\frac 12 ||B||.
 \]
 Thus we have
 \begin{equation}\label{use}
 \int_{B(R_2)\backslash B(R_1)}H\,d\Sigma\geq \frac 12 \eps_0 (R_2-R_1)
 \end{equation}
 provided that $R_2>R_1$ are large enough. 
 
We will create suitable test functions using the techniques similar to ~\cite{dek,cek-1,ll-1,ll-2,ll-4}.
 Let $\phi\in C_0^\infty(\Sigma\backslash B(\frac R2))$ be a smooth function such that
 \[
 \phi\equiv 1\quad \text{on } B(2R)\backslash B(R),\quad
 \int_\Sigma|\nabla\phi|^2 d\Sigma<\eps_1,
 \]
 where  $\eps_1\to 0$ as $R\to\infty$. The existence of such a function $\phi$ is guaranteed by the parabolicity of $\Sigma$. Then we have, as in Remark~\ref{rem21}, that
 \[
 Q(\phi\chi,\phi\chi)<4a\eps_1+2a\pi^2\int_{\Sigma\backslash B(R/2)} K\phi^2d\Sigma.
 \]
 Since $K$ is integrable, for any $\eps_2>0$, there is an $R_0>0$ such that if $R>R_0$, we have
 \[
 Q(\phi\chi,\phi\chi)<\eps_2.
 \]

 Now let's consider a function $j\in C_0^\infty (B(\frac 53 R)\backslash B(\frac 43 R))$. Consider the function $j\chi(t)t$, where $j$ is a smooth function on $\Sigma$ such that $j\equiv 1$ on $B(\frac{19}{12}R)\backslash B(\frac{17}{12}R)$; $0\leq j\leq 1$; and $|\nabla j|<2$. Then there is an absolute  constant $C_1$, such that
 \[
 Q(j\chi(t)t,j\chi(t)t)\leq C_1 R^2.
 \]
Finally, let's consider $Q(\phi\chi(t),j\chi(t)t)$. Since ${\rm\,\, supp}j\subset \{\phi\equiv 1\}$, by~\eqref{q1}, $Q_1(\phi\chi(t),j\chi(t)t)=0$. Let
\[
\sigma=-\int^a_{-a} \chi'(t)\chi(t)tdt>0.
 \]
 Then
 \[
 Q(\phi\chi(t),j\chi(t)t)=-\sigma\int_\Sigma j d\Sigma.
 \]
 Let $\eps>0$. Then we have
 \[
 Q(\phi\chi(t)+\eps j\chi(t)t, \phi\chi(t)+\eps j\chi(t)t)
 <\eps_2-2\eps\sigma\int_\Sigma jd\Sigma+\eps^2C_1R^2.
\]
By ~\eqref{use}, we have
\[
Q(\phi\chi(t)+\eps j\chi(t)t, \phi\chi(t)+\eps j\chi(t)t)
 <\eps_2-\frac 13\eps\sigma R+\eps^2C_1R^2.
 \]
 If 
 \[
 \eps_2<\frac{\sigma^2}{36 C_1},
 \]
 then there is a suitable $\eps>0$ such that 
 \[
 Q(\phi\chi(t)+\eps j\chi(t)t, \phi\chi(t)+\eps j\chi(t)t)
 <0.
 \]
 Thus $\sigma_0<\kappa^2$.
 
 \qed

 \section{Further Discussions.}
 We proved the following more general
 
 \begin{theorem}
We assume that $\Sigma$ satisfies
\begin{enumerate}
\item
The isopermetric inequality holds. That is, there is a constant $\delta_1>0$ such that if $D$ is a domain  in $\Sigma$, we have
\[
(length(\pa D))^2\geq \delta_1 Area(D).
\]
\item There is another positive constant $\delta_2>0$ such that for any compact set $K$ of $\Sigma$, there is a curve $C$ outside the set $K$ such that if $\vec \gamma$ is one of its normal vector in $\Sigma$, then there is a vector $\vec a$ such that
\[
\langle \vec \gamma,\vec a\rangle\geq\delta_2>0
\]
for some fixed vector $\vec a\in R^3$.
\end{enumerate}
Then the ground state exists.
\end{theorem}

 {\bf Proof.}
 We let $\phi$ be a smooth function such that ${\rm supp}\,\phi\subset B(R)\backslash B(r)$ for $R>>R/4>>4r>>r>0$ large. We also assume that on $B(R/2)\backslash B(2r)$, $\phi\equiv 1$. Let $\eps_0>0$ be a positive number to be determined later such that
 \[
 \int_\Sigma|\nabla\phi|^2\leq \eps_0,\qquad
 \int_\Sigma|K|\phi^2\leq \eps_0.
 \]
 Note that $\eps_0$ is independent of $R$.
 
 We let $\chi=\cos\frac{\pi}{2a} t$. Then there is a constant $C$ such that
 \[
 Q(\phi\chi,\phi\chi)<C\eps_0.
 \]
 
 Let $C$ be a curve outside the compact set $B(4r)$ satisfying the condition in the theorem. We let $R$ big enough that $C\subset B(R/4)$. 
 
 In order to construct the test functions, we let $\rho$ be the cut-off function
 such that $\rho=1$ if $t\leq 0$ and $\rho=0$ if $t\geq 1$ and we assume that $\rho$ is decreasing. Near the curve $C$, any point $p$ has a coordinate $(t,s)$, where $s\in C$ from the exponential map. 
 To be more precise, let $(x_1,x_2)$ be the the local coordinates near $C$ such that locally $C$ can be represented by $x_1=0$. Let the Riemannian metric under this coordinate system be
\[
g_{11} (dx_1)^2+2g_{12} dx_1 dx_2+g_{22} (dx_2)^2.
\]
The fact that $\vec\gamma$ is a normal vector implies that if
\[
\vec\gamma=\gamma_1\frac{\pa}{\pa x_1}+\gamma_2\frac{\pa}{\pa x_2},
\]
then
\[
\gamma_1g_{12}+\gamma_2g_{22}=0.
\]
Let $\sigma_t(x_2)$ be the geodesic lines starting from $x_2\in C$ with initial vector $\vec \gamma$. Then $\sigma_t$ is the exponential map. The Jacobian of the map at $t=0$ is
\[
\begin{pmatrix}
\gamma_1&\gamma_2\\
0&1
\end{pmatrix}
\]
In particular, $\gamma_1\neq 0$ since the map must be nonsingular. A simple computation shows that $\nabla t=\gamma_1g^{1j}\frac{\pa}{\pa x_i}$. Thus $\nabla t$ is proportional to $\vec \gamma$.

 Let $\phi_1$ be a cut-off function   such that $\phi_1\equiv 1$ on $B(R/4)\backslash B(4r)$ and ${\rm supp}|\,(\phi_1)\subset B(R/2)\backslash B(2r)$.

We define $\tilde \rho(p)=\phi_1\rho(t/\eps_1)$, where $\eps_1$ is a positive constant to be determined. WLOG, let $\vec a$ be the $z$-direction in the Euclidean space.

Let $\vec n$ be the normal vector of $\Sigma$. Let $n_z$ be the $z$-component of $\vec n$.  We compute the following term
 $Q(\phi\chi, \tilde\rho n_z\chi_1)$, where $\chi_1=t\cos\frac{\pi}{2a}t$. First
 $Q_1(\phi\chi, \tilde\rho n_z\chi_1)=0$ becasue ${\rm supp}(\tilde\rho n_z)$ is contained in the area where $\phi\equiv 1$. On the other hand, since $\chi\chi_1$ is an odd function, we have
 \[
 Q_2(\phi\chi, \tilde\rho n_z\chi_1)=
 -\int_\Sigma H\phi\tilde\rho n_z d\Sigma\int_{-a}^a(\chi'\chi_1' t-\kappa^2\chi\chi_1 t) dt.
 \]
 A straight computation shows that
 \[
 C_1=\int_{-a}^a(\chi'\chi_1' t-\kappa^2\chi\chi_1 t) dt=-1/2\neq 0.
 \]
 Furthermore, we have $Hn_z=\Delta z$. As a result, we have
 \[
 -\int_\Sigma H\phi\tilde\rho n_z d\Sigma=\int_\Sigma\nabla z\nabla\tilde\rho=\int_{\{t\leq\eps_1\}}\nabla z\nabla\tilde\rho
 \]
 (Note that $\phi\equiv 1$ on the points we are interested).
 We have the following Taylar expansion:
 
 \[
 \nabla z\nabla\tilde\rho(t,x_2)= \nabla z\nabla\tilde\rho(0,x_2)+O(t)
 \]
 Since $\int _{\{t\leq\eps_1\}}O(t)/\eps_1=O(\eps_1)\,Length(C)$, ,we have
 \[
 \int_\Sigma\nabla z\nabla\tilde\rho\geq(\delta_2-O(\eps_1))\,Length(C)
 \]
 We choose $\eps_1$ small enough, then we have
 \[
 \int_\Sigma\nabla z\nabla\tilde\rho\geq\frac 12\delta_2\,Length(C)
 \]
 
  If we let $\eps\rightarrow 0$, then that above becomes
 \[
 -\int_\Sigma H\phi\rho n_z d\Sigma=\int_\Sigma\nabla z\nabla\rho
 \geq \delta_2\,Length (C).
 \]
 
 Finally, we have $|\rho n_z|+|\nabla(\rho n_z)|\leq 2$, thus we have
 \[
 Q(\rho n_z\chi_1,\rho n_z\chi_1)\leq C Area(D),
 \]
 where $D$ is the domain $C$ enclosed. To summary, for any $\eps<0$, we have
 \[
 Q(\phi\chi+\eps\rho n_z\chi_1,\phi\chi+\eps\rho n_z\chi_1)
 \leq C\eps_0+2\eps C_1\delta_2\, Length(C)+C\eps^2 Area(D).
 \]
 Using the isopermetric inequality, we know that if $\eps_0<\delta_1\delta_2^2/C^2$ is small enough, then
 \[
Q(\phi\chi+\eps\rho n_z\chi_1,\phi\chi+\eps\rho n_z\chi_1)
<0
\]
which proves the theorem. 
 
 \qed

Using the 
 same  proof, we can prove the following:
 
 \begin{theorem} Using the same notations as in Conjecture~\ref{ccon}, we assume further that
 \[
 ||B||(x)\leq C/dist (x,x_0),
 \]
 where $x_0\in\Sigma$ is a reference point of $\Sigma$. Then Conjecture~\ref{ccon} is true.
 \end{theorem}
 
\qed

\bibliographystyle{abbrv}   
\bibliography{pub,unp,2007}   
\end{document}